\documentclass[12pt]{article}
\usepackage{amsmath,amssymb,amsfonts}

\newcommand{\R}{\mathbb{{R}}}
\newcommand{\Z}{\mathbb{{Z}}}

\newcommand{\lo}{\longrightarrow}

\newcommand{\ov}{\overline}

\newcommand{\para}{\paragraph}
\begin{document}
\begin{center}
{\bf\large  Lusternik-Schnirelmann Theory for a Morse Decomposition}\\[1cm]
M.R. Razvan
\\ \noindent {\it \small Institute for Studies in Theoretical Physics and
Mathematics\\ \small P.O.Box: $19395-5746$, Tehran, IRAN}\\
E-Mail:
 razvan@karun.ipm.ac.ir
\end{center}

\vspace*{5mm}
\begin{abstract}
Let $\varphi^t$ be a continuous flow on a metric space $X$ and
$I$ be an isolated invariant set with an index pair $(N,L)$ and a
Morse decomposition $\{M_i\}^n_{i=1}$. For every category $\nu$
on $N/L$, we prove that $\nu(N/L)\leq \nu([L])+\sum_{i=1}^n
\nu(M_i)$. As a result if $\varphi^t|_I$ is gradient-like and $X$
is semi-locally contractible, then $\varphi^t$ has at least
$\nu_H(h(I))-1$ rest points in $I$
 where $h(I)$ is the Conley index of $I$
and $\nu_H$ is the Homotopy Lusternik-Schnirelmann category.
\end{abstract}



\section{Introduction}
Conley's homotopy index theory was first introduced as a
generalization of Morse theory \cite{C} and it was indeed a
landmark in this field and related subjects. From the view point
of critical points, Morse theory concerns the relations between
the topology of a manifold and the number of analytically distinct
critical points of smooth functions on it. However the natural
problem is the number of geometrically distinct critical points
and this is investigated by Lusternik-Schnirelmann theory. An
important result in this theory  is that if there is a smooth
function $f:M\lo \R$ on a compact boundaryless manifold $M$ with
$n$ critical points, then $M$ is covered by $n$ contractible open
subsets \cite{B}. The above result can be generalized for compact
metric spaces as follows: If a compact locally
 contractible metric space $X$ admits a gradient-like flow
with $n$ rest points, then $X$ is covered by $n$ contractible
open subsets of $X$.

The compactness assumption is crucial in the above results.
Indeed every noncompact manifold admits a smooth function without
critical points. But in some applications of critical point
theory, we deal with a compact subset of a noncompact space
\cite{CZ1}. In \cite{C} and \cite{CZ2}, the results of Morse
theory were generalized for every isolated invariant set. In this
paper we use Conley index theory to obtain Lusternik-Schnirelamnn
results for isolated invariant sets.

We need some basic results from Conley index theory which are
presented in the next section. Then in section 3, we define the
concept of category and prove a Ljusternik-Schnirelmann  result
in Conley index theory. In section 4, we consider a well-known
example of Lusternik-Schnirelmann category and we obtain the
results of \cite{P}.

\section{Conley Index Theory}
Let $\varphi^t$ be a continuous flow on a metric space $X$. An
isolated invariant set is a subset $I$ of $X$ which is the maximal
invariant set in a {\bf compact} neighborhood of itself. Such a
neighborhood is called an isolating neighborhood.

\para{Definition.} A Morse decomposition of $I$ is a collection
$\{M_i\}^n_{i=1}$ where each $M_i$ is an isolated invariant subset
of $I$ and for all $x\in I-\bigcup_{i=1}^n M_i$ there exist
$i,j\in \{1,\cdots , n\}$ such that $i<j$, $\alpha (x) \in M_i$
and $\omega (x)\in M_j$.

Let $V$ be any isolating neighborhood for $I$. We define
\begin{center}
$I^+=I^+(V)=\{x\in V \ |\  \varphi^{[0,\infty)} (x) \subset N\}$,
\\ $I^+_j=I^+(N)=\{x\in N\ |\  \varphi^{[0,\infty)} (x) \subset N \
\text{and} \ \omega (x)\subset M_1\cup \cdots M_j\}$
\end{center}
for $j=1,\cdots , n$. In \cite{CZ2} it is proved that $I_j$ is
compact and for every $x\in N$ with $\varphi^{[0,\infty)}
(x)\subset N$, there is a $j\in\{1,\cdots , n\}$
 such that
$\omega (x) \subset M_j$ (See \cite{S}, Lemma 3.8).

In order to define the Conley  index of an isolated invariant set, we follow
\cite{RS}.
 Given a compact pair $(N,L)$ we define the induced semi-flow
on $N/L$ by: $$\varphi_{\#}^t: N/L \lo N/L ,  \ \ \varphi_{\#}^t
(x)=\left\{ \begin{array}{ll} \varphi^t (x) & \text{if} \
\varphi^{[0,t]} (x) \subset N-L\\ \left[ L \right] &
\text{otherwise.}
\end{array} \right. $$
In \cite{RS} it is proved that $\varphi^t_{\#}$ is continuous if and only if:

i) $L$ is positively invariant relative to $N$, i.e.
$$x\in L, \varphi^{[0,t]} (x) \subset N \Rightarrow \varphi^{[0,t]} (x) \subset
L.$$

ii) Every orbit which exits $N$ goes through $L$ first:
$$x\in N,\varphi^{[0,\infty)} (x) \not\subset N \Rightarrow \exists_{t\geq 0} \
\text{with} \ \varphi^{[0,t]} (x) \subset N, \varphi^t (x) \in L,$$
or equivalently if $x\in N-L$ then there is a $t>0$ such that $\varphi^{[0,t]}
(x) \subset N$.

 An index pair for an isolated invariant set $I\subseteq X$ is
a compact pair $(N,L)$ in $X$ such that $\overline{N-L}$ is an
isolating neighborhood for $I$ and the semi-flow $\varphi^t_{\#}$
induced by $\varphi^t$ is continuous.In \cite{C}, \cite{S}
 and \cite{RS} it is shown that every  isolated invariant set $I$ admits
an index pair $(N,L)$ and the homotopy type of the pointed space
$N/L$ is independent of the choice of the index pair. The Conley
index of $I$ is the homotopy type of $N/L$.

\para{Note.} We shall not distinguish between $N-L$ and
$N/L-\{[L]\}$.

\para{Lemma 2.1.} If $A$ is a compact subset of $N/L$ with
$A\cap I^+ =\varnothing$, then there exist a $T\in \R ^+$ such
that $\varphi^T_{\#} (A)=[L]$. In particular $[L]$ admits a
contractible neighborhood in $N/L$. \\ {\bf Proof.} Since $I$ is
the maximal invariant set in $\ov{N-L}$, for every $x\in A$,
there is a $t\in \R ^+$ such that $\varphi^t(x)\not\in \ov{N-L}$.
Thus $\varphi^t(y)\not\in \ov{N-L}$ for every $y$
 in a neighborhood of $x$. Now
by compactness,
 we can find a $T\in \R ^+$ with $\varphi^{[0,T]}(x)\not\subset \ov{N-L}$
for every $x\in A$ and the proof is complete. It also shows that
$(\varphi_{\#}^t)^{-1} ([L])$ is a compact contractible
neighborhood of $[L]$ for sufficiently large amounts of $t$.
$\square$

\para{Lemma 2.2.} Let $I$ be an isolated invariant set with an
isolating neighborhood $N$ and a Morse decomposition
$\{M_i\}^n_{i=1}$. If $A$ is a closed subset of $I_j^+-I_{j-1}^+$
and $U$ is a neighborhood of $M_j$, then there exists a $T\in \R
^+$ such that $\varphi^t(A)\subset U$, for every $t\geq T$.
\\ {\bf Proof.} We may assume that $U$ is a compact neighborhood
of $M_i$ such that $U\subset N$ and $U\cap I_{j-1}^+=\varnothing$.
Now for every $x\in U\cap (I_j^+ -\stackrel{\circ}{U})$ there is a
$t\in \R ^+$ such that $\varphi^{-t} (x)\not\in U$. Since $U\cap
(I_j^+-\stackrel{\circ}{U})$ is compact, there exists a $T\in \R
^+$ such that $\varphi^{[-T,0]} (x) \not\subset U$ for every
$x\in U\cap (I_j^+-\stackrel{\circ}{U})$. Now choose a
neighborhood $V$ of $M_j$ such that $\varphi^{[0,T]}(V)\subset
\stackrel{\circ}{U}$. We claim that $\varphi^{[0,\infty)} (V\cap
I_j^+) \subset \stackrel{\circ}{U}$. Suppose the contrary, then
$\varphi^t (x)\not\in\stackrel{\circ}{U} $ for some $x\in V\cap
I_j^+$ and $t\in \R ^+$. Since $\varphi^0(x)=x\in
\stackrel{\circ}{U}
 \cap I_j^+$, there is a $t_0<t$
such that $\varphi^{t_0} (x)\in U \cap
(I_j^+-\stackrel{\circ}{U})$ for the first time, i.e.
$\varphi^s(x)\in \stackrel{\circ}{U}$ for $s\in [0,t_0)$. Since
$\varphi^{[0,T]}(V) \subset \stackrel{\circ}{U}$, we have
$T<t_0<t$. Thus $\varphi^{t_0-T} (x)\not\in U$ which is a
contradiction with $\varphi^s(x)\in \stackrel{\circ}{U}$ for
$0<s<t_0$. Now for every $x\in A$, there is a $t\in \R ^+$ such
that $\varphi^t (x) \in V$ and since $A$ is compact, there is a
$T\in \R ^+$ such that $\varphi^s(x) \in V$ for some $s\in [0,T]$.
Since $\varphi^s(x)\in V\cap I_j^+$ and $\varphi^{[0,\infty)}
(V\cap I_j^+)\subset U$, we conclude that $\varphi^t(A)\subset
U$, for every $t\geq T$. $\square$

\section{Lusternik-Schnirelmann Theory}
Let $M$ be a topological space. A category on $M$ is a map $\nu:
2^M\lo \Z \cup \{+\infty\}$ which satisfies the following
axioms:\\ i) If $A\subset B$, then $\nu(A)\leq \nu(B)$.\\ ii) $\nu
(A\cup B)\leq \nu (A)+\nu (B)$.\\ iii) For every $A\subset M$,
there is an open set $U\subset M$ with $A\subset U$ and
$\nu(A)=\nu(U)$.\\ iv) If $f:M\lo M$ is homotopic to the identity
$id_M$, then $\nu(A)\leq \nu (f(A))$ for every $A\subset M$.

It has been shown that if $\nu$ is a category on a compact metric
space  $X$  satisfying the following axiom: \\ v) If A consists of
a single point, then $\nu (A)=1$,\\ then every gradient-like flow
on $X$ possesses at least $\nu(X)$ lest points (cf.\cite{MS}). The
following theorem gives a generalization of this result and also
the resuls of \cite{P} and \cite{Ru}.

\para{Theorem 3.1.} Let $\varphi^t$ be a continuous flow on a
metric space $X$ and $I$ be an isolated invariant set for
$\varphi^t$ with a Morse decomposition $\{M_i\}^n_{i=1}$. If
$(N,L)$ is an index pair for $I$ and $\nu$ is a category on
$N/L$, then $\nu(N/L)\leq \nu ([L])+\sum_{i=1}^n \nu (M_i)$.
Moreover, if $L= \varnothing$, then $\nu(N)\leq \sum_{i=1}^n \nu
(M_i)$. \\ {\bf Proof.} Recall that for every $1\leq j\leq n$,
$I_j^+$ is compact where $I_j^+=\{x\in N| \varphi^{[0,\infty)}
(x)\subset N$ and $\omega (x)\subset M_1 \cup \cdots \cup M_j\}$,
$I^+_j\subseteq I_{j+1}^+$ and $I^+_n=I^+=\{x\in
N|\varphi^{[0,\infty)} (x) \subset N\}$. We construct $V_j
\supset M_j$ open  in $N/L$ such that $\nu(V_j)=\nu(M_j)$ and
$I_j\subset \bigcup_{i=1}^j V_i$ for $1\leq j \leq n$ by
induction. For $j=1$, by axioms (i) and (iii), there is an open
set $U_1\subset N-L$ such that $M_1\subset U_1$ and $\nu
(M_1)=\nu (U_1)$. Now by Lemma 2.2, there is a $t_1 \in \R ^+$
such that $\varphi^{t_1}(I^+_1)\subset U_1$ and hence
$I^+_1\subset V_1:=\{x\in N|\varphi^{[0,t_1]} (x)\subset N \
\text{and} \ \varphi^{t_1} (x) \in U_1\}.$ It is easy to see that
$V_1$ is an open set in $N$, $M_1\subset V_1\subset N-L$ and
$M_1\subset \varphi_{\#}^{t_1} (V_1) \subset U_1$. Since $
\varphi_{\#}^{t}$ is homotopic to identity, by axioms (i) and
(iv) $\nu(M_1)\leq \nu (V_1)\leq \nu (\varphi_{\#}^t (V_1))\leq
\nu (V_1)=\nu (M_1)$ and hence $\nu (V_1)=\nu (M_1)$. $V_1$ is
the desired open set in $N/L$.

Now suppose that we have constructed $V_1,\cdots , V_{j-1}$ open
in $N/L$ with $M_i\subset V_i$ and $I^+_{j-1}\subset
\bigcup_{i=1}^{j-1} V_i$. Similar to the case $j=1$, we choose an
open set $U_j\subset N-L$ with $\nu(U_j)=\nu(M_j)$. Now
$I^+_j-\bigcup_{i=1}^{j-1} V_i$ is a compact subset of $I_j^+
-I_{j-1}^+$ and by Lemma 2.1 there exists a $t_j \in \R ^+$ such
that $\varphi^{t_j}(I^+_j-\bigcup_{i=1}^{j-1} V_i)\subset U_j$
and hence $$I^+_j-\bigcup_{i=1}^{j-1} V_i\subseteq V_j:=\{ x\in
N|\varphi^{[0,t_j]} (x) \subset N \ \text{and} \ \varphi^{t_j}
(x)\in U_j\}$$ It is not hard to see that  $V_j$ is open in $N$,
$M_j\subset V_j \subset N-L$ and $M_j\subset
 \varphi_{\#}^{t_j} (V_j)\subset U_j$. Thus $\nu (M_j)=\nu (V_j)$
and $I^+_j\subseteq \bigcup_{i=1}^{j} V_i$.

Now if $L=\varnothing$, then $N=I^+=I^+_n=\bigcup_{i=1}^n V_i$
and by axiom (ii), $\nu(N)\leq \sum_{i=1}^n \nu
(V_i)=\sum_{i=1}^n \nu (M_i)$. If $L\neq \varnothing$, then by
Lemma 2.1. there is a $T\in \R ^+$ such that $\varphi^T_{\#}
(N/L-\bigcup_{i=1}^n V_i)=[L]$. This shows that $\nu
(N/L-\bigcup_{i=1}^n V_i)=\nu ([L])$ by axiom (iv) and now by
axiom (ii) $$\nu (N/L)\leq \nu ([L])+\sum_{i=1}^n \nu
(V_i)=\nu([L])+\sum_{i=1}^n \nu (M_i). \ \square$$

\para{Corollary 3.2.} Let $I$ be a compact invariant subset of $X$
and $\{M_i\}^n_{i=1}$ be a Morse decomposition for $I$. Then for
every category $\nu$ on $X$, $\nu(I)\leq \sum_{i=1}^n \nu(M_i)$
\\ {\bf Proof.} Consider the flow $\varphi^t|_I$ on compact
metric space $I$. Then $(I,\varnothing)$ is an index pair for $I$
with respect to $\varphi^t|_I$. Now the above argument shows that
there exist $V_1,\cdots , V_n$ open in $I$ and $t_1,\cdots ,
t_n\in \R ^+$ such that $\bigcup_{i=1}^n V_i=I$ and $\varphi^{t_i}
(V_i)\subset U_i$ where $U_i$ is an open subset of $I$ with
$\nu(M_i)=\nu(U_i)$. Since $\varphi^{t_i}$ is defined on $X$ and
homotopic to $id_X$, we have $\nu(V_i)\leq \nu(U_i)=\nu(M_i)$ by
axiom (iv) and hence $\nu(I)\leq \sum_{i=1}^n \nu(V_i)\leq
\sum_{i=1}^n \nu(M_i)$. $\square$

Now suppose that $X$ is a compact manifold and $\varphi ^t$ is the
gradient flow of a smooth function with $n$ critical points. Then
in the above corollary, each $U_i$ can be chosen a disk and it
follows that $X$ can be covered by $n$ open disks. Moreover if
$f:X\lo \R$ is a smooth function with $n$ critical value
$c_1<\cdots <c_n$ and $M_i$ is the set of critical points in
$f^{-1}(c_i)$, then we can use the above corollary to obtain the
results of \cite{Ru}. Similarly if $X$ is a locally contractible
metric space and $\varphi ^t$ is a gradient-like flow with $n$
rest points, then $X$ can be covered by $n$ contractible open
sets.

\para{Theorem 3.3.} Let $\varphi^t$ be a continuous flow and
$I$ be an isolated invariant set for $\varphi^t$ such that
$\varphi^t|_I$ is gradient-like. Then for every index pair $(N,L)$
for $I$ and every category $\nu$ on $N/L$ which satisfies axiom
(v) in $I$, $\varphi^t$ has at least $\nu (N/L)-\nu([L])$
 rest points in $I$.\\
{\bf Proof.} We may assume that $\varphi^t|_I$ has a finite number
of rest points $x_1,\cdots , x_n$. Since $\varphi^t|_I$ is
gradient-like, these rest point give a Morse decomposition of $I$.
Now by Theorem 3.1., $\nu (N/L)\leq \nu ([L])+\sum_{i=1}^n \nu
(x_i)$. Since $\nu$ satisfies axiom (v) in $I$,
 $\nu (N/L)\leq \nu([L])+n$. Thus the number of
rest points of $\varphi^t|_I$ is not less than
$\nu(N/L)-\nu([L])$. $\square$

\section{ HLS Category}

Let $M$ be a topological space. A subset $A\subset M$ is called
contractible in $M$ if the inclusion map $A\lo M$ is homotopic to
a constant. We define Homotopy Lusternik-Schnirelmann category as
follows:

\para{Definition.}
 The HLS-category $\nu_H(A)=\nu_H(A,M)$ of a subset $A\subset M$ is
defined to be the minimum number of open sets contractible in $M$ required to
cover $M$. If such a cover does not exist, we set $\nu_H(A)=\infty$ and if it
exists, $A$ is called $H$-categorizable (in $M$).

A subset $A\subset M$ is H-categorizable if and only if $\nu_H(\{x\})=1$
  for every $x\in A$. Thus
 $\nu_H$ satisfies axiom (v) if and only if $M$ is $H$-categorizable.
It is easy to see that $\nu_H$ satisfies axioms (i)-(iii). The
following result \cite{J} gives a generalization of axiom (iv).

\para{Lemma 4.1.} If $Y$ dominates $X$, i.e. there are
continuous maps
 $f:X\lo Y$ and
 $g:Y\lo X$ with  $g\circ f\sim
id_X$), then for every $H$-categorizable subset $A\subset Y$,
$f^{-1} (A)$ is $H$-categorizable and $\nu_H(f^{-1}(A))\leq
 \nu_H(A)$. In particular if $Y$ is $H$-categorizable, then
so is $X$ and $\nu_H(X)\leq \nu_H(Y)$. \\
{\bf Proof:} It is enough to prove that for every open set
$U\subset Y$ contractible in $Y$, $f^{-1}(U)$ is contractible in
$X$. Consider the following commutative diagram in which $i$ and
$j$ are inclusion maps:

\unitlength=1mm \special{em:linewidth 0.4pt} \linethickness{0.4pt}
\begin{picture}(83.33,43.33)
\put(50.00,20.00){\vector(1,0){25.00}}
\put(50.00,40.00){\vector(1,0){25.00}}
\put(75.00,35.00){\vector(-1,0){25.00}}
\put(45.00,22.00){\vector(0,1){8.00}}
\put(80.00,22.67){\vector(0,1){8.00}}
\put(80.00,37.00){\makebox(0,0)[cc]{$Y$}}
\put(62.67,43.33){\makebox(0,0)[cc]{$f$}}
\put(62.67,31.33){\makebox(0,0)[cc]{$g$}}
\put(45.00,37.00){\makebox(0,0)[cc]{$X$}}
\put(40.00,26.33){\makebox(0,0)[cc]{$i$}}
\put(40.00,20.00){\makebox(0,0)[cc]{$V=f^{-1}(U)$}}
\put(62.67,16.33){\makebox(0,0)[cc]{$f|_V$}}
\put(80.00,20.00){\makebox(0,0)[cc]{$U$}}
\put(83.33,26.67){\makebox(0,0)[cc]{$j$}}
\end{picture}

\vspace*{-1cm}

$f\circ i=j\circ f|_V \Rightarrow g \circ f \circ i = g\circ j
\circ f|_V \Rightarrow i \sim g\circ j \circ f|_V \sim$ constant.
$\square$

\para{Corollary 4.2.} Let $X$ and $Y$ be topological spaces of
the same homotopy type. If $X$ is $H$-categorizable then so is
$Y$ and
$\nu_H(X)=\nu_H(Y)$.\\
{\bf Proof:} We have two maps $f:X\lo Y$ and $g:Y\lo X$ such that
$f\circ g \sim id_Y$ and $g\circ f \sim id_X$. The above lemma
shows that $Y$ is $H$-categorizable
  and $\nu_H(Y)\leq \nu_H(X)$. Now since $Y$ is $H$-categorizable,
then so is $X$
 and $\nu_H(X)\leq \nu_H(Y)$. $\square$

\para{Remark 4.3.} Let $I$ be an isolated set with a Morse
decomposition $\{M_i\}^n_{i=1}$ and two index pairs $(N,L)$ and
$(N',L')$. Then there is a flow-defined homotopy equivalence
between $N/L$ and $N'/L'$ which leaves each $M_i$ invariant
\cite{S}. Thus by Lemma 4.1., $\nu_H (M_i, h(I))$, $\nu_H
(I,h(I))$ and $\nu_H (h(I),h(I))$ make sense.  Moreover by
 Lemma 2.1, $[L]$ admits a contractible neighborhood in $N/L$
and hence $\nu_H ([L], N/L)=1$.

\para{Theorem 4.4.} Let $I$ be an isolated invariant set with a
Morse decomposition $\{M_i\}^n_{i=1}$. If
$\nu_H(M_i,h(I))<\infty$ for $1\leq i \leq n$, then $h(I)$ is
$H$-categorizable and $\nu_H (h(I), h(I))\leq 1 +\sum_{i=1}^n
\nu_H (M_i, h(I))$.

The above theorem is an immediate consequence of Theorem 3.1. In
order to apply theorem 3.3, we need $h(I)$ to be
$H$-categorizable. We introduce a class of metric spaces in which
the Conley  index of every isolated invariant set is
$H$-categorizable.

\para{Definition.}
 A topological space $X$ is called semi-locally contractible if for every
$x\in X$ and open set $U\subseteq X$ with $x\in U$, there exists a
neighborhood $V$ of $x$ such that $x\in V\subset U$ and $V$ is
contractible in $U$.

\para{Lemma 4.5.} Let $I$ be an isolated invariant set in a
semi-locally contractible metric space $X$. Then $\nu_H(I,
h(I))<\infty$
 and hence $h(I)$ is
$H$-categorizable. \\ {\bf Proof.} Suppose that $(N,L)$ is an
index pair for $I$. Since $X$ is semi-locally contractible, every
$x\in I$ admits an open set contractible in $N-L$. Since $I$ is
compact, $\nu_H(I, h(I))<\infty$
 and by Theorem 4.4, case $n=1$, $h(I)$ is $H$-categorizable. $\square$

\para{Theorem 4.6.} Let $I$ be an isolated invariant set for a
continuous flow $\varphi^t$ on a semi-locally contractible metric
space $X$. If $\varphi^t|_I$ is gradient-like, then
$\varphi^t$ has at least $\nu_H (h(I))-1$ rest points in $I$.\\
{\bf Proof.} $\nu_H$ satisfies axiom (v) on $h(I)$ by lemma 4.6.
and we can use Theorem 3.4. $\square$

\para{Acknowledgment.} The author would like to thank Institute for
studies in theoretical Physics and Mathematics, IPM ,  for
supporting this research.

\end{document}